\documentclass{amsart}
\usepackage{amssymb}
\usepackage{amsmath, amsfonts}
\usepackage{setspace, csquotes}
\usepackage[super]{nth}
\usepackage{xcolor}
\usepackage{mathtools}
\usepackage[colorinlistoftodos]{todonotes}

\newtheorem{theorem}{Theorem}[section]

\newtheorem{lemma}[theorem]{Lemma}

\theoremstyle{definition}
\newtheorem{definition}[theorem]{Definition}

\numberwithin{equation}{section}

\begin{document}

\title[Moment of double exponential sums]
{Moment of double exponential sums}
\author{Nilanjan Bag}
 \address{Department of Mathematics, Thapar Institute of Engineering and Technology, Patiala, Punjab 147004, India}
\email{nilanjanb2011@gmail.com}
\author{Dwaipayan Mazumder}
 \address{Chennai Mathematical Institute, Kelambakkam, Siruseri, Tamil Nadu 603103, India}
\email{dwmaz.1993@gmail.com}
\date{November 10, 2025}
\begin{abstract}
This paper is devoted to finding moments of double exponential sums with monomials over arbitrary sets and intervals in finite fields. The study of such sums dates back to the work of Heath-Brown, who studied such sums in a work on least square-free numbers in an arithmetic progression.  
\end{abstract}
\maketitle
\section{Introduction and statements of the results}
Let $p$ be an odd prime. Consider the field of $p$ elements which in standard convention denoted as $\mathbb{F}_p$. Take two subsets $\mathcal{M}$ and $\mathcal{N}$ in $\mathbb{F}_p$. In this paper, we focus our study on 
\begin{align}\label{mainsum}
\mathcal{S}(\mathcal{M},\mathcal{N},r,s;k)=\sum_{x\in\mathcal{M}}\left|\sum_{y\in\mathcal{N}}{\bf{e}}_p\left(ax^ry^s\right)\right|^{2k},
\end{align}
 where ${\bf{e}}_p(z)=exp(\frac{2\pi i z}{p})$. This study complements the first author and Shaprlinski's 2021 work \cite{Bag-Igor-1} on double exponential sums. Such sums are important in number theory for its connections in finding smallest square free numbers in arithmetic progression, in the study of divisor function and also its contribution to Brune-Titchmarsh theorem. There are many varieties of results. For example for $k=1/2$ and $r=s=-1$, from a result of Bourgain \cite{Bour-1}, one can see that the above sum can achieve a non-trivial bound $MNp^{-\delta}$, provided $MN\geq p^{1/2+\varepsilon}$, for some $\varepsilon>0.$ Beyond the square root thresh hold it can achieve the same bound for $M\geq p^{1/8}$ and $N\geq p^{5/12+\varepsilon}$. For more details, see \cite[Theorem 7]{Bour-2}.
 \par A subset $\mathcal{A}\subset \mathbb{F}_p^{\times}$ is called an interval if it is of the form
 \begin{align*}
 \mathcal{A}=\{A+1,...,A+N\},
\end{align*} 
where when $A=0$, we call $\mathcal{A}$ as {\it initial interval}, otherwise it is called an {\it arbitrary interval}. For $s=-2$ and $k=1/2$ with any fixed integer $r$, one can see that 
\begin{align}\label{bound-1}
\mathcal{S}(\mathcal{M},\mathcal{N},r,-2;1/2)\ll MN\left(\frac{p}{MN^{2\ell/(1+\ell)}}\right)^{1/{2\ell}}p^{o(1)}
\end{align}
for any fixed integer $\ell\geq 1$, an arbitrary set $\mathcal{M}$ with cardinality $M$ and an initial interval $\mathcal{N}$ with cardinality $N\leq p^{(\ell+1)/2\ell}$. For more details see \cite[eq 21, Lemma 1]{heath-brown}. In \cite{Bag-Igor-1} the first author and Shparlinski proposed a new approach to bound such sums which applies in more generality and works for any arbitrary interval $\mathcal{N}$, instead of only initial interval. For complex weights ${\boldsymbol{\alpha}}=(\alpha_x)_x$, a weighted version of double exponential sums can be taken as
\begin{align*}
    \mathcal{S}_{\alpha}(\mathcal{M},\mathcal{N},r,s)=\sum_{x\in\mathcal{M}}\alpha_x\sum_{y\in\mathcal{Y}}{\bf{e}}_p\left(ax^ry^s\right).
\end{align*}
Trivially we have 
\begin{align*}
    | \mathcal{S}_{\alpha}(\mathcal{M},\mathcal{N},r,s)|\leq ||\boldsymbol{\alpha}||_{\infty}|\mathcal{S}(\mathcal{M},\mathcal{N},r,s;1/2)|.
\end{align*}
As $\mathcal{M}$ is an arbitrary set, the power $r$ becomes irrelevant and one can consider $r=1$. 
In \cite{Bag-Igor-1}, authors prove the following bound
\begin{align*}
\mathcal{S}_{\alpha}(\mathcal{M},\mathcal{N},1,s)\leq ||\alpha||_1^{1-\frac{1}{\ell}}||\alpha||_{\infty}^{1/\ell}M^{1/{2\ell}}N\left(N^{-1/\ell}p^{(\ell+1)/{2\ell^2}}+1\right)p^{o(1)},
\end{align*}
which when compared with \eqref{bound-1} is weaker but applies to more generality over arbitrary interval.
\par This paper is devoted to generalizing the result of Bag and Shparlinski in \cite{Bag-Igor-1} to even power moments \eqref{mainsum}. 
We are interested in the case when $s=-1$ and $r=1$ and $2k$-power moment of double exponential sum. We denote this sum as $S(\mathcal{M}, \mathcal{N}, k)$ and is defined as 
\begin{align*}
    \mathcal{S}(\mathcal{M},\mathcal{N}, k)=\sum_{m\in\mathcal{M}}\left|\sum_{n\in\mathcal{N}}{\bf{e}}_p\left(a\overline{m}n\right)\right|^{2k},
\end{align*}
where $\mathcal{M}\subseteq \mathbb{F}_p^{\times}$. First  fix the notation $$\underbrace{ \Sigma + \Sigma + \cdots +\Sigma }_{k \text{ times}} := k\Sigma, ~\text{and}~ \Sigma := \mathcal{N - N}.$$ Then  for $\mathcal{M}$ an arbitrary interval and $\mathcal{N}$ an arbitrary subset, we have the following theorem.

\begin{theorem}\label{MT1}
Let $\mathcal{M}\subseteq\mathbb{F}_p^{\times}$ be an arbitrary interval of length $M$ and $\mathcal{N}\subset\mathbb{F}_p^{\times}$ be an arbitrary subset of cardinality $N$ with $\Sigma\neq \mathbb{F}_p$. Then for any fixed integer $r\geq 2$, we have 
\begin{align*}
\mathcal{S}(\mathcal{M},\mathcal{N},k)\ll p^{1/2r^2}M^{1-1/r}p^{(k-2)/2}|\Sigma|^{k/2}|k\Sigma|^{1-1/2r}p^{o(1)},
\end{align*}
which is non-trivial for $N> p^{\frac{1}{2}-\frac{1}{4rk}(1-\frac{1}{r})}.$
 \end{theorem}
\section{Results from Additive combinatorics}
In this section we  define Fourier bias and state a lemma which helps us in the proof of our main theorem. We need Fourier analysis on the finite additive group $Z$. We start with the following definition.
\begin{definition}
A bilinear form on an additive group $Z$ is a map $(\xi, x) \longmapsto \xi\cdot x$ from $Z \times Z$ to $\mathbb{R}/\mathbb{Z}$, which is a homomorphism in each of the variables $\xi, x$ separately.
\end{definition}
A bilinear form is non-degenerate if for every non-zero $\xi$ the map $x \longmapsto \xi\cdot x$ is not identically zero, and similarly for every non-zero $x$ the map $\xi \longmapsto \xi \cdot x$ is not identically zero. We also say the form is symmetric if $\xi \cdot x = x \cdot \xi$. \par
For $A\subseteq Z$ (which we call an additive set), $\mathbb{P}_Z(A)$ denote the density or probability of $A$  and is defined as  
$$
\mathbb{P}_Z(A): = \frac{|A|}{|Z|}.
$$
Next we define additive energy as,
\begin{definition}
    If $A$ and $B$ are additive sets with ambient group $Z$, we define the additive energy $E(A, B)$ between $A$ and $B$ to be the quantity
$$
E(A , B) := |\{(a , a' , b , b') \in A \times A \times B \times B : a + b = a' + b'\}|.
$$
\end{definition}

Let $\mathbb{C}^Z$ denote the space of all complex-valued functions $f:Z\longrightarrow \mathbb{C}$. If $f\in \mathbb{C}^Z$, we define the mean or expectation of $f$ to be the quantity 
$$
\mathbb{E}_Z(f) = \mathbb{E}_{x\in Z}f(x) := \frac{1}{|Z|}\sum_{x\in Z}f(x).
$$
For $A\subseteq Z,$ the function $1_A(\xi)$ is the characteristic function on $A$. If $f \in \mathbb{C}^Z$, the Fourier transform $\hat{f} \in \mathbb{C}^Z$ is expressed by the formula 
\begin{equation*}
\hat{f}(\xi) := \mathbb{E}_{x \in Z}f(x)\overline{{e(x.y)}}.
\end{equation*}
where $e(x)=exp(2 \pi i x).$
\subsection{Fourier Bias}
\par We now define Fourier bias of a subset $A\subseteq Z\setminus\{0\}$.
\begin{definition}
    Let $Z$ be a finite additive group. For $A\subseteq Z\setminus\{0\}$, the Fourier bias ${||A||}_u$ of the set $A$ is defined to be the quantity
$$
||A||_u := \sup\limits_{\xi\in Z\setminus{\{0\}}} |\hat{1}_A(\xi)|
$$
\end{definition}
We note that $||A||_u = 0$ if and only if $A = Z$.
The following lemma describes the connection between Fourier bias and sum sets.
\begin{lemma}[\textbf{Uniformity implies large sum sets}]\cite[Lemma 4.13]{tao}\label{Bias}
    Let $n\geq 3$, and let $A_1, A_2,..., A_n$ be additive sets in a finite additive group $Z$. Then for any $x \in Z$ we have 
\begin{align*}
\left|\frac{1}{|Z|^{n-1}}\left|\{(a_1,...,a_n)\in A_1\times\cdots \times A_n:x=a_1+\cdots +a_n\}\right|-P_{Z}(A_1)\cdots P_{Z}(A_n)\right|\\
\leq ||A_1||_u\cdots ||A_{n-2}||_uP_{Z}(A_{n-1})^{1/2}P_{Z}(A_{n})^{1/2}.
\end{align*} 
\end{lemma}
Let $A$ be an additive set in a finite additive group $Z$. In \cite[equation 4.23]{tao}, it is claimed that
\begin{align}\label{bound2}
    ||A||_u^4 \leq \frac{1}{|Z|^3}E(A , A) - \mathbb{P}_Z(A)^4 \leq ||A||_u^2\mathbb{P}_Z(A).
\end{align}
Here we give a proof. We note for $A \subset Z$,
\begin{align*}
&\sum\limits_{\xi\in Z }\left|\frac{1}{|Z|}\sum_{a \in A}e(\xi\cdot a)\right|^4 \\=& \frac{1}{|Z|^4}\sum\limits_{\xi\in Z }\left(E(A , A) + \sum\limits_{\substack{a_1. a_2, a_3. a_4 \in A \\ a_1 + a_3 \neq a_2 + a_4}}e(\xi\cdot ((a_1 + a_3) - (a_2 + a_4)))\right)\\
=& \frac{1}{|Z|^3}E(A, A).
\end{align*}
Also 
$$
\sum\limits_{\xi\in Z }\left|\frac{1}{|Z|}\sum_{a \in A}e(\xi\cdot a)\right|^4 = \frac{|A|^4}{|Z|^4} + \sum\limits_{\xi\in Z}\left|\frac{1}{|Z|}\sum_{a \in A}e(\xi\cdot a)\right|^4 \geq \frac{|A|^4}{|Z|^4} + ||A||_u^4.
$$
So altogether we get
\begin{align}\label{lowerbound}
||A||_u^4 \leq \frac{1}{|Z|^3}E(A , A) - \mathbb{P}_Z(A)^4.
\end{align}
If we invoke Lemma \ref{Bias} with $n = 4, A_1 = A_3 = A, A_2 = A_4 = -A$ and $x = 0$ we get 
\begin{align}\label{uupperbound}
\frac{1}{|Z|^3}E(A , A) - \mathbb{P}_Z(A)^4 \leq ||A||_u^2\mathbb{P}_Z(A).
\end{align}
Relations \eqref{lowerbound} and \eqref{uupperbound} implies that 
$$
||A||_u^4 \leq \frac{1}{|Z|^3}E(A , A) - \mathbb{P}_Z(A)^4 \leq ||A||_u^2\mathbb{P}_Z(A),
$$
This completes the proof of the \eqref{bound2}. Now an immediate consequence of \eqref{bound2} is 
\begin{align}\label{upperbound}
||A||_u \leq {\mathbb{P}_Z(A)}^{1/2}.
\end{align}

The inequalities \eqref{upperbound} plays a crucial role in the proof of our main theorem. In our context we take $Z = \mathbb{Z}/p\mathbb{Z}$ and one can show that $x \cdot \xi := x\xi/p$ is symmetric and non-degenerate bilinear form on $\mathbb{Z}/p\mathbb{Z}$.

\section{Proof of Theorem \ref{MT1}}
We begin with breaking the modulus power as
\begin{align}\label{bound-final}
    \mathcal{S}(\mathcal{M}, \mathcal{N}, k) &= \sum\limits_{m \in \mathcal{M}}|\sum\limits_{n \in\mathcal{N}}{\bf{e}}_p(\overline{m}n)|^{2k}\notag\\
    & = \sum\limits_{m \in \mathcal{M}}\quad\sum\limits_{n_1, n_2, ...n_{2k} \in \mathcal{N}} {\bf{e}}_p(\overline{m}(\sum_{i = 1}^{k}n_i - \sum_{j = r+1}^{2k} n_j))\notag\\
    & \ll \max\limits_{0\leq l\leq p-1}\rho(l)\sum\limits_{n \in \mathfrak{N}}|\sum\limits_{m \in \mathcal{M}}{\bf{e}}_p(\overline{m}n)|,
\end{align}
where 
\begin{align*}
\rho(l) = \sum\limits_{\substack{n_1, n_2,....n_{2k}\in \mathcal{N}\\\sum_{i = 1}^{k} n_i - \sum_{j = r+1}^{2k} n_j \equiv l \bmod p}} 1
\end{align*}
and
$$
\mathfrak{N} = \underbrace{ \Sigma + \Sigma + \cdots +\Sigma }_{k \text{ times}} := k\Sigma, \quad \Sigma = \mathcal{N - N}.
$$
\subsection{Estimation of $\rho(l)$:}
\par To get an upper bound for $\rho(l)$ we use Lemma \ref{Bias}. Taking $A_i = \Sigma$ and $Z = \mathbb{Z}/p\mathbb{Z}$, we get 
$$
\left|\frac{\rho(l)}{|Z|^{k-1}} - \frac{|\Sigma|^k}{|Z|^k}\right| \leq {||\Sigma||_u}^{(k-2)}\cdot \frac{|\Sigma|}{|Z|}.
$$
Hence $\rho(l)\leq |Z|^{k-1}\left\{||\Sigma||_u^{k-2}\frac{|\Sigma|}{|Z|}+\frac{|\Sigma|^k}{|Z|^k}\right\}=|Z|^{k-2}|\Sigma|||\Sigma||^{k-2}_u+\frac{|\Sigma|^k}{|Z|}.
$
Using \eqref{upperbound} in the above, we deduce
\begin{align}\label{numberofsol}
    \rho(l)&\leq |Z|^{(k-2)/2}|\Sigma|^{k/2}+\frac{|Z|^k}{|Z|}\leq 
     2p^{(k-2)/2}|\Sigma|^{k/2},
\end{align}
as $|\Sigma|\leq |Z|.$ 
\subsection{Estimation of $\sum\limits_{n \in \mathfrak{N}}|\sum\limits_{m \in \mathcal{M}}{\bf{e}}_p(\overline{m}n)|$:}
\par Now to treat the inner sum we use "ab-shifting" technique. We rename the inner sum as 
\begin{align*}
    \mathfrak{G}=\sum\limits_{n \in \mathfrak{N}}|\sum\limits_{m \in \mathcal{M}}{\bf{e}}_p(\overline{m}n)|
\end{align*}
Let $\mathcal{M}$ be the arbitrary interval of the form 
\begin{align*}
    \mathcal{M}=\{I+1,...I+M\}.
\end{align*}
We fix two positive parameters $A$ and $B$ such that $AB\leq M$. Consider two sets $\mathcal{A}$ and $\mathcal{B}$, which are set of integers in $[A/2,A]$ and $[B/2,B]$. Then we follow the same strategy as followed in \cite[Section 4.1]{KMS}, to get
\begin{align}\label{bound3}
 \mathfrak{G}\ll \frac{\log p}{AB}\sum_{n\in\mathfrak{N}}\sum_{a\in\mathcal{A}}\sum_{m\in \mathcal{M}' }\left|\sum_{b\in\mathcal{B}}\eta(b){\bf{e}}_p\left((\overline{m\overline{a}+b})\overline{a}n\right)\right|,
\end{align}
for some complex numbers $\eta(b)$ with $|\eta(b)|=1$ and $\mathcal{M}'=[I-M,I+2M]$. For $s$ and $t$ in $\mathbb{F}_p^{\times}$, let $\nu(s,t)$ be the cardinality of the set 
\begin{align*}
    \{(a,m,n)\in\mathcal{A}\times \mathcal{M}'\times \mathfrak{N}|\overline{a}m=s, \overline{a}n=t \}.
\end{align*}
Hence rewriting \eqref{bound3}, we have 
\begin{align*}
    \mathfrak{G}\ll \frac{\log p}{AB}\sum_{t\in\mathbb{F}_p^{\times}}\sum_{s\in\mathbb{F}_p^{\times}}\nu(s,t)\left|\sum_{b\in\mathcal{B}}\eta(b){\bf{e}}_p\left((\overline{s+b})t\right)\right|.
\end{align*}
Write $\nu(s,t)=\nu(s,t)^{(r-1)/r}\nu(s,t)^{1/r}$, and then using H\"{o}lder inequality with weights $r/(r-1)$, $2r$ and $2r$, we get 
\begin{align}\label{bound-G}
    \mathfrak{G}^{2r}\leq \frac{1}{(AB)^{2r}}\mathcal{R}_1^{2r-2}\mathcal{R}_2\mathcal{S}p^{\varepsilon},
\end{align}
for some $\varepsilon>0$.
Here \begin{align*}
    &\mathcal{R}_1=\sum_{s,t\in\mathbb{F}_p^{\times}}\nu(s,t)\\
    &\mathcal{R}_2=\sum_{s,t\in\mathbb{F}_p^{\times}}\nu(s,t)^2\\
    &\mathcal{S}=\sum_{s,t\in\mathbb{F}_p^{\times}}\left|\sum_{b\in\mathcal{B}}\eta(b){\bf{e}}_p\left((\overline{s+b})t\right)\right|^{2r}.
\end{align*}
Trivially 
\begin{align}\label{bound4}
    \mathcal{R}_1\ll AM |\mathfrak{N}|.
\end{align}
Now from \cite[p.334]{FI}, we have 
\begin{align}\label{bound5}
     \mathcal{R}_2\ll AM |\mathfrak{N}|p^{\varepsilon}.
\end{align}
To estimate $\mathcal{S}$, we first open up the inner sum and then following same technique as \cite[p.234]{Bag-Igor-1} to get 
\begin{align}\label{bound6}
    S\ll pB^r+B^{2r}.
\end{align}
We omit the details for reasons of brevity. Put $A=Mp^{-1/r} $ and $B=p^{1/r}$. Then comibing \eqref{bound4}, \eqref{bound5}, \eqref{bound6} in \eqref{bound-G}, we get
\begin{align}\label{bound8}
   \mathfrak{G}&\ll \frac{p^{\epsilon/{r}} }{AB}\left(AM|\mathfrak{N}|\right)^{1-\frac{1}{2r}}B\notag\\
&=
     \frac{p^{\epsilon/{r}+1/r} }{M}\left(M^2|\mathfrak{N}|p^{-1/r}\right)^{1-\frac{1}{2r}}\notag\\
   & \ll p^{\epsilon/{r}+1/2r^2}M^{1-1/r}|k\Sigma|^{1-1/2r}.
\end{align}
Putting bounds \eqref{numberofsol} and \eqref{bound8} in \eqref{bound-final}, we finally have
\begin{align*}
    \mathcal{S}(\mathcal{M}, \mathcal{N}, k)  &\ll  p^{\epsilon/{r}+1/2r^2}M^{1-1/r}p^{(k-2)/2}|\Sigma|^{k/2}|k\Sigma|^{1-1/2r}.
\end{align*}
Now comparing for 
\begin{align*}
   p^{(k-2)/2}|\Sigma|^{k/2}|k\Sigma|^{1-1/2r}< N^{2k},
\end{align*}
we get $N> p^{\frac{1}{2}-\frac{1}{4rk}(1-\frac{1}{r})}.$ This completes the proof of Theorem \ref{MT1}.
\section{Acknowledgement}
The first author would like to thank Department of Mathematics, Thapar Institute of Engineering and Technology and the second author would like to thank Chennai Mathematical Institute for providing excellent working conditions. 
\par During the preparation of this article, D.M. was supported by the National Board of Higher Mathematics post-doctoral fellowship (No.: 0204/10/(8)/2023/R\&D-II/2778).



\begin{thebibliography}{99}
\bibitem{Bag-Igor-1}
 N. Bag, I. E. Shparlinski, {\it Bounds of some double exponential sums},  Journal of Number Theory 219 (2021),
228-236.

\bibitem{Bour-1}
 J. Bourgain, {\it More on the sum-product phenomenon in prime fields and its applications},
 Int. J. Number Theory, 1(1) (2005),
1-32.
\bibitem{Bour-2}
J Bourgain, MZ Garaev, {\it 
Sumsets of reciprocals in prime fields and multilinear Kloosterman sums},
Izv. Ross. Akad. Nauk, Ser. Mat. 78 (2014) 9–72 (in Russian); translation in Russ. Acad. Sci. Izv.
Math. 78 (2014) 656–707.

\bibitem{FI}
 J. B. Friedlander and H. Iwaniec, {\it Incomplete Kloosterman Sums and a Divisor Problem}, Annals of Mathematics, Vol. 121, No. 2 (Mar., 1985), pp. 319-344.
\bibitem{heath-brown}
D. R. Heath Brown, {\it The least square-free number in an arithmetic progression},
 J. Reine Angew.
Math. 332 (1982) 204–220.
\bibitem{KMS}
 E. Kowalski, P. Michel, W. Sawin, {\it Stratification and averaging for exponential sums: bilinear forms
with generalized Kloosterman sums}, Ann. Sc. Norm. Pisa 21 (2020) 1453–1530.
\bibitem{tao}
T. Tao, V. H. Vu, {\it Additive combinatorics},
Combridge studies in advanced mathematics, 105.

\end{thebibliography}
\end{document}